\documentclass[letterpaper, 10 pt, conference]{ieeeconf}  

\IEEEoverridecommandlockouts                              

\usepackage[maxfloats=256]{morefloats}
\usepackage{mathtools}
\usepackage{amsmath,amssymb,amsfonts,color}
\usepackage{subfigure}
\usepackage{color}
\usepackage{float}

\newcommand{\cA}{\mathcal{A}}
\newcommand{\cC}{\mathcal{C}}

\newcommand{\R}{\mathbb{R}}

\newcommand{\cL}{\mathcal{L}}

\newcommand{\cJ}{\mathcal{J}}

\newcommand{\bU}{\mathbf{U}}

\newcommand{\bW}{\mathbf{W}}
\newcommand{\bb}{\mathbf{b}}

\newcommand{\bo}{\mathbf{0}}

\newcommand{\bA}{\mathbf{A}}
\newcommand{\bQ}{\mathbf{Q}}
\newcommand{\bB}{\mathbf{B}}
\newcommand{\bR}{\mathbf{R}}
\newcommand{\bK}{\mathbf{K}}

\newcommand{\bv}{\mathbf{v}}

\newcommand{\bx}{\mathbf{x}}
\newcommand{\by}{\mathbf{y}}
\newcommand{\bu}{\mathbf{u}}
\newcommand{\bff}{\mathbf{f}}

\newtheorem{proposition}{Proposition}

\title{\LARGE \bf
Gradient-augmented Supervised Learning of Optimal Feedback Laws Using State-dependent Riccati Equations 
}

\author{Giacomo Albi$^{1}$, Sara Bicego$^{1}$ and Dante Kalise$^{2}$
\thanks{DK was supported by a public grant as part of the Investissement d’avenir project,reference ANR-11-LABX-0056-LMH, LabEx LMH, and by the UK Engineering and Physical Sciences Research Council (EPSRC) grants EP/V04771X/1, EP/T024429/1, and EP/V025899/1.}
\thanks{$^{1}$Giacomo Albi and Sara Bicego are with the Department of Computer Science, University of Verona, Strada le Grazie 15 - 37134 Verona, Italy
        {\tt\small giacomo.albi@univr.it, sara.bicego@studenti.univr.it}}%
\thanks{$^{2}$Dante Kalise is with the School of Mathematical Sciences, University of Nottingham, University Park, Nottingham NG7 2RD, United Kingdom 
        {\tt\small dante.kalise@nottingham.ac.uk}}%
}

\begin{document}

\maketitle
\thispagestyle{empty}
\pagestyle{empty}

\begin{abstract}
A supervised learning approach for the solution of large-scale nonlinear stabilization problems is presented. A stabilizing feedback law is trained from a dataset generated from State-dependent Riccati Equation solves. The training phase is enriched by the use gradient information in the loss function, which is weighted through the use of hyperparameters. High-dimensional nonlinear stabilization tests demonstrate that real-time sequential large-scale Algebraic Riccati Equation solves can be substituted by a suitably trained feedforward neural network.
\end{abstract}

\section{INTRODUCTION}
A large class of control problems in fluid flow control, consensus dynamics, and power networks, among many others, can be cast a optimal stabilization problems sharing two distinctive features in the dynamics: nonlinearity, and a high-dimensional state space. The natural control-theoretical framework to address these problems is via optimal stabilization using dynamic programming and Hamilton-Jacobi-Bellman (HJB) partial differential equations (PDEs). Unfortunately, the HJB PDE arising in nonlinear control is a first-order fully nonlinear equation with no general explicit solution. Moreover, the overwhelming computational complexity associated to the solution of high-dimensional HJB PDEs poses a formidable challenge limiting the applicability of traditional numerical methods such as finite differences or finite elements to very low-dimensional control systems.

The numerical approximation high-dimensional HJB PDEs arising in deterministic optimal control is a topic that has been the subject of extensive research. Without attempting an exhaustive literature review on the topic, effective  computational approaches to this problem include the use of sparse grids~\cite{axelsparse}, tree structure algorithms~\cite{allatree}, max-plus methods ~\cite{maxplus2}, polynomial approximation ~\cite{poli1,KK17} and tensor decomposition methods~\cite{tensor1,tensor2,tensor3,tensor4,tensor5}. These grid-based schemes are complemented with recent works making use of artificial neural networks~\cite{ml1,ml2,ml3,ml4,ml5}.

In this paper, we propose a computational method for the solution of large--scale optimal stabilization problems for nonlinear dynamics avoiding the solution of the HJB PDE through a supervised learning approach. This idea dates back to \cite{BTB00}, where the synthesis of feedback controls by interpolating finite horizon open-loop solves was proposed. More recently, this problem has been studied in \cite{kang1,kang2} using a sparse grid interpolant, in \cite{nakazim,nakazim2} using deep neural networks, and in \cite{AKK21} through sparse polynomial regression. Similarly, the works \cite{chow1,chow3,chow4} make use of representation formulas for HJB PDEs along with fast convex optimization solvers.

Many of the aforementioned works exploit the relation between the Hamilton-Jacobi-Bellman PDE and necessary optimality conditions through Pontryagin's Maximum Principle (PMP) in finite horizon control. Under convexity and smoothness assumptions, the PMP system represents the characteristic curves of the HJB PDE, and the value function of the problem can be computed at a given space-time point by solving a two-point boundary value problem. Unfortunately, such an interpretation is not readily available for infinite horizon optimal control, which is the case of interest for asymptotic stabilization of nonlinear dynamics. 
 
 The methodology proposed in the present work circumvents the direct solution of the HJB PDE and the lack of PMP-like representation formula for the value function by resorting to State-dependent Riccati Equations (SDRE) \cite{BLT07,C97}. In the SDRE framework, after casting the nonlinear dynamics in semilinear form, a feedback control is obtained by a sequential solution of Algebraic Riccati Equations (ARE) along the trajectory. Under certain stabilizability conditions, this feedback law generates a locally asymptotically stable closed-loop and approximates the optimal feedback law from the HJB PDE. However, the main computational bottleneck of the SDRE approach is the availability of a sufficiently fast ARE solver to be called at an arbitrarily high rate. In this paper, we propose a supervised learning approach to replace the real-time solution of AREs by the use of an artificial neural network for the feedback law generated by the SDRE approach. We demonstrate that, through an adequate choice of network architectures, and including the use of gradient information of the model in the training, is it possible to accurately recover the SDRE feedback law for high-dimensional nonlinear stabilization problems.
 
 The rest of the paper is organized as follows. In Section II we describe the nonlinear optimal stabilization setting, and in Section III we discuss its solution via the SDRE approach. In Section IV we discuss its numerical approximation through supervised learning, to conclude in Section V with a computational assessment for two nonlinear, high-dimensional tests, presenting concluding remarks in Section VI.

\section{Infinite Horizon Optimal Feedback Control}
We study the design of feedback laws for asymptotic stabilization through infinite horizon optimal control:
\begin{equation}\label{ocp}
	\underset{\bu(\cdot)\in \bU}{\min}\cJ(\bu(\cdot),\bx_0):=\int\limits_0^\infty \bx^{\top}\!(s)\bQ\bx(s)\,+\,\bu^{\top}\!(s)\bR\bu(s)\,ds\,,
\end{equation}
subject to nonlinear, control-affine dynamics of the form
\begin{align}\label{dynamics}
	\dot \bx(t)=\bff(\bx(t))+\bB(\bx(t))\bu(t)\,,\quad\bx(0)=\bx_0\,,
\end{align}
where $\bx(t)=(x_1(t),\ldots,x_n(t))^{\top}\in\R^n$ denotes the state of the system, $\bu(\cdot)\in\bU=\{\bu(t):\, \R^+\rightarrow \R^m, \text{measurable}\}$ is an unbounded control variable, $\bQ\in\R^{n\times n}$ is a symmetric positive semidefinite matrix, and  $\bR\in\R^{m\times m}$ is symmetric positive definite. The system dynamics $\bff(\bx):\R^n\rightarrow\R^n$ and the control operator $\bB(\bx):\R^n\rightarrow\R^{n\times m}$ are assumed to be $\cC^1(\R^n)$ and, without loss of generality, such that $\bff(\bo)=\bo$ and $\bB(\bo)=\bo$. The optimal feedback law for the control problem \eqref{ocp} is synthesized using Dynamic Programming. For this, we define the value function of the control problem
\begin{equation}
	V(\bx)=\underset{\bu(\cdot)\in \bU}{\inf}\cJ(\bu(\cdot),\bx)\,,
\end{equation}
which in turn satisfies a first-order, static, nonlinear Hamilton-Jacobi-Bellman PDE
\begin{align}\label{hjb}
	\nabla V(\bx)^{\top}\!\bff(\bx)&-\frac14\nabla V(\bx)^{\top}\!\bW(\bx)\nabla V(\bx)\!+\!\bx^{\top}\!\bQ\bx=0\,,
\end{align}
where $\bW(\bx)=\bB(\bx)\bR^{-1}\bB(\bx)^{\top}$. After solving for $V(\bx)$, the optimal feedback is given by
\begin{equation}\label{feedback}
	\bu(\bx)=-\frac{1}{2}\bR^{-1} \bB(\bx)^{\top} \nabla V(\bx)\,.
\end{equation}
The main difficulty when applying the dynamic programming approach to optimal feedback synthesis resides in the solution of the HJB PDE \eqref{hjb}. This is a nonlinear PDE cast in the state-space of the system dynamics, with a dimension that can be arbitrarily high. Perhaps the most successful instance of a solution to this problem is the \emph{linear quadratic regulator } (LQR), where, under the additional assumption that the free dynamics are linear, $\bff(\bx)=\bA\bx$ and $\bB(\bx)=\bB$, and making the ansatz $V(\bx)=\bx^{\top}\!\Pi\bx$ with $\Pi\in\R^{n\times n}$ leads to
\begin{equation}\label{feedbackare}
	\bu(\bx)=-\bK\bx=-\bR^{-1} \bB^{\top} \Pi\bx\,,
\end{equation}
where $\Pi$ is a positive definite solution of the Algebraic Riccati Equation (ARE)
\begin{equation}\label{are}
	\bA^{\top}\Pi+\Pi \bA-\Pi \bB\bR^{-1}\bB^{\top}\Pi+\bQ=0\,.
\end{equation}
There are different methods which utilize the solution of the ARE above to generate a sub-optimal feedback control for local stabilization of nonlinear dynamics. Most notably, solving \eqref{are} with $(\bA,\bB(\bo))$ where $\bA_{ij}=\frac{\partial \bff_i(\bx)}{\partial x_j}|_{\bx=0}$ leads to a linear feedback operator $\bK_0$ which can effectively stabilize states in a vicinity of the origin. In the following, we discuss the synthesis of nonlinear feedback control laws by a sequential solution of AREs.

\section{State-Dependent Riccati Equation}\label{secsdre}
Having a representation of the nonlinear dynamics in semilinear form
\begin{equation}\label{sdress}
	\dot \bx=\bA(\bx)\bx+\bB(\bx)\bu(t)\,,
\end{equation} 
we approximate the synthesis of the optimal feedback control following the State-dependent Riccati Equation (SDRE) approach. Formally, the solution of the nonlinear optimal control problem \eqref{ocp} is associated to an ARE where the operators are state-dependent
\begin{align}\label{sdre}
	\bA^{\top}\!(\bx)\Pi(\bx)&\!+\!\Pi(\bx) \bA(\bx)\!-\!\Pi(\bx) \bW(\bx)\Pi(\bx)\!+\!\bQ\!=\!0\,,
\end{align}
and analogously, the feedback \eqref{feedbackare} is also expressed through a state-dependent gain operator $\bK(\bx)$
\begin{equation}\label{sdref}
	\bu(\bx)=-\bK(\bx)\bx=-\bR^{-1}\bB^{\top}(\bx)\Pi(\bx)\bx\,.
\end{equation}
Aiming at directly solving \eqref{sdre} for a general high-dimensional operator $\Pi(\bx)$ leads to the same difficulties already present in \eqref{hjb}. Instead, we assume the operator $\Pi(\bx)$ is a positive definite matrix in $\R^{n\times n}$, so that for a fixed $\bx$, solving \eqref{sdre} effectively reduces to problem to an ARE. We can benefit from this SDRE framework by applying it in a receding horizon fashion. Given a current state $\bar\bx$ along a trajectory, we solve \eqref{sdre} for $\Pi(\bar\bx)$ by freezing every operator accordingly, recovering the feedback $\bu(\bx)=-\bK(\bar\bx)\bx$, to then evolve the controlled dynamics for a reduced time frame, after which we update the state of the system and recompute the feedback law. This approach leads to two natural questions: establishing conditions under which the SDRE approach generates an asymptotically stable closed-loop, and the design of effective computational methods for the fast solution of SDREs of potentially large scale. Regarding the first question, we recall the following proposition on asymptotic stability of the closed-loop generated by the SDRE approach \cite{BLT07}. 
\begin{proposition}Assume a nonlinear system
	\begin{equation}
		\dot \bx(t)=\bff(\bx(t))+\bB(\bx(t))\bu(t)\,,
	\end{equation}
	where $\bff(\bx)$ is $\cC^1$ for $\|\bx\|\leq\delta$, and $\bB(\bx)$ is continuous. If $\bff(x)$ is parametrized in the form $\bff(\bx)=\bA(\bx)\bx$, and the pair $(\bA(\bx),\bB(\bx))$ is stabilizable for every $\bx$ in a non-empty neighbourhood of the origin $\Omega\subset\mathcal{B}_{\delta}(\bo)$, then the closed-loop dynamics generated by the feedback law \eqref{sdref} are locally asymptotically stable.
\end{proposition}
Assuming the stabilizality hypothesis above, the main bottleneck in the implementation of the SDRE approach is the availability of an ARE solver sufficiently fast for real-time feedback control. Here, we assume an ARE solver is readily available, however, it is not suitable for real-time control. In order to circumvent this difficulty, we follow a supervised learning approach, as we explain in the following section.
\section{Gradient-Augmented Supervised Learning for Optimal Feedback Laws}
The SDRE \eqref{sdre} is solved offline for a set of training states, denoted by $\mathcal{X}_t$, which is used for training a suitable artificial neural network (ANN) which is then implemented for real-time control. The use of ANNs for SDREs has been explored for learning the matrix-valued operator $\Pi(\bx)$ in \eqref{sdre}, see e.g. \cite{WANG98}. Here, we propose two alternatives:

\paragraph{\textbf{Learning} $\bu(\bx)$}: we train a model for the vector-valued feedback law $\bu(\bx): \R^n\to\R^m$ upon a set of $N_s$ training states $\mathcal{X}_t:=\{\bx^{(i)}\}_{i=1}^{Ns}$, the solution of the corresponding $\Pi(\bx)$, and the controls $\bu(\bx)$ via \eqref{sdref}.
\paragraph{\textbf{Learning} $V(\bx)$} we train a model for the scalar function $V(\bx): \R^n\to\R$ from $V(\bx)=\bx^{\top}\Pi(\bx)\bx$ and its gradient $\nabla V(\bx)=2\Pi(\bx)$, where $\Pi(\bx)$ is a positive definite solution of \eqref{sdre} for each $\bx\in\mathcal{X}_t$. The feedback law is then expressed as $	\bu(\bx)=-\frac{1}{2}\bR^{-1} \bB(\bx)^{\top} \nabla V(\bx)$.

Both alternatives are a direct supervised learning formulation of the SDRE approach, with the sole objective of synthesizing a feedback requiring a reduced number of operations for online implementation. However, the second approach links the solution of the SDRE with finding a function $V(\bx)$ which approximates the solution of the original HJB equation \eqref{hjb}. As discussed in \cite{Astolfi2020}, there is a direct equivalence between HJB, SDRE, and ARE in the linear-quadratic case. For the general nonlinear case, the ansatz $V(\bx)=\bx^{\top}\Pi(\bx)\bx$ with $\Pi(\bx)$ generated from the SDRE approximates the solution of the HJB PDE only in neighbourhood of the origin. However, this idea is instrumental from a computational viewpoint. The advantage of the second formulation resides in the training of a scalar function, for which both function and gradient values are available. This shall be reflected in the choice of gradient-augmented loss functions for training.

\subsubsection{Network architecture}
The approximation task is carried out using feedforward neural networks (FNN), with information flowing from the input nodes to the output without generating any cycles or loops. FNNs approximate a function $f(\cdot)$ by a chain of compositions 
\begin{align}
    f(\bx) &\approx f_{\theta}(\bx) = l_{M} \circ ... \circ l_{2} \circ l_{1}(\bx),
\end{align} 
where each layer $l_{m}$ is defined as $l_{m}(\by) = \sigma_m(\bA_m \by + \bb_m)$, $\bA_m$ are the weight matrices, $\bb_m$ are the bias vectors and $\sigma_m(\cdot)$ is a nonlinear \emph{activation function} applied component-wise. Standard choices for $\sigma(\cdot)$ are the ReLU function $\sigma(x) = max(0,x)$ and $\sigma(x) = tanh(x)$. The activation function in the hidden layers needs to be chosen accordingly with the valuation of the model's goodness of fit, and the last layer is typically assumed to be linear, thus $\sigma_{M}(x) = x$. 

Considering a data set $\mathcal{T} = \{\bx^{(i)}, f(\bx^{(i)})\}_{i=1}^{N_s}$, the NN is trained over the parameters $\theta = \{\bA_m, \bb_m\}_{m=1}^M$ to best approximate the target $f(\bx)$, i.e. minimizing the \emph{loss} between the approximation $f_{\theta}(\bx^{(i)})$ of the model and the true values $f(\bx^{(i)})$ for every $\bx^{(i)} \in \mathcal{T}$:
\begin{equation}\label{loss}
    \min_{\theta} \; \cL(f(\bx), f_{\theta}(\bx))
\end{equation}
where the loss function $\cL$ evaluates how well $f_{\theta}(\cdot)$ models the given dataset $\mathcal{T}$.  The goodness of fit of the trained NN $f_{\theta}$, within a set $\mathcal{T}' = \{\bx^{(j)}\}_{j=1}^{N_v}$ can be measured by the \emph{coefficient of determination} 
\begin{equation}
	r^{2} = 1 - \frac{\sum_{j=1}^{N_v} \|f(\bx^{(j)}) - f_{\theta}(\bx^{(j)})\|^2}{\sum_{j=1}^{N_v}  \|f(\bx^{(j)}) - \bar{f}\|^2},
\end{equation}
where $\bar{f} = \frac{1}{n} \sum_{j=1}^{N_v}  f(\bx^{(j)})$. This coefficient typically ranges in $[0,1]$;  an $r^2$ of $1$ indicates that the model approximations perfectly fit the data, while values below $0$ suggest the trained model to fit the data worse than a horizontal hyperplane.

We search for an approximation of the feedback control $\bu(\bx)$, for which we consider two different approaches: to build a model $\bu_{\theta}(\cdot)$ having $\bu(\cdot)$ itself as target variable, $\bu(x) \approx \bu_{\theta}(x)$, or to describe it through a FNN $V_{\theta}$ approximating $V(\cdot)$, on top of which we add a feedback layer
\begin{equation}\label{u_dV}
    \bu(x) \approx \bu_{V}(x) = -\frac{R^{-1}B^T\nabla V_{\theta}(x)}{2}.
\end{equation}
An accurate approximation of $\nabla V(\cdot)$ is essential for calculating a reasonable $\bu_{V}(x)$. Here we deal with this through automatic differentiation, which allow us to compute exact gradients of $V_{\theta}$ in an efficient way. In this case, our training is not limited to pointwise valuations of $V(\bx)$, but also includes the discrepancy between the true gradient $\nabla V(\bx)$ and its approximation $\nabla V_{\theta}$. This is done choosing an ad hoc loss function $V_{\theta}$. 

\subsubsection{Loss function}
The training of the neural network for $u_{\theta}$ is done through a standard loss function: the mean squared error (MSE)
\begin{equation}
    \cL_{0}(\bu, \bu_{\theta}) := \frac{1}{N_s}\sum_{i=1}^{N_s}\|u(\bx^{(i)}) - u_{\theta}(\bx^{(i)})\|^2,
\end{equation}
averaging the squared difference between approximation and actual observations. 

For the training of $V_{\theta}$, we consider instead  
\begin{equation}\label{loss1}
    \cL_{1}(V, V_{\theta}) =  \mu_V  \cL_{0}(V, V_{\theta}) + \mu_{dV}  \cL_{0}(\nabla V, \nabla V_{\theta})\,.
\end{equation}
This loss function represents a compromise between the fitting functional $\cL_{0}(V, V_{\theta})$ and the gradient regulation $\cL_{0}(\nabla V, \nabla V_{\theta})$, suitably weighted thanks to $\mu_V$ and $\mu_{dV}$.

\section{Numerical Experiments}
We assess the neural network approximation for feedback laws in two different tests. The control laws to be approximated rely on the pointwise solution of the SDRE \eqref{sdre}, for which we resort to the \texttt{lqr} routine in MATLAB. The samples for training were generated by solving (\ref{ocp})-(\ref{dynamics}) for initial condition vectors $\mathcal{X}_t = \{\bx^{(i)}\}_{i = 1}^{N_s}\in \Omega\subset\mathbb{R}^n$, being populated using Halton quasi-random sequences in $[0,1]^n$. 

Once the solution of the SDRE is computed for each sample $\bx^{(i)} \in \mathcal{X}_t$, the training set $\{\bx^{(i)},\bu(\bx^{(i)})\}_{i=1}^{N_s}$ for $u_{\theta}$  can be computed as in (\ref{sdref}), while the ANN $V_{\theta}$ is trained upon an enriched dataset, containing both the value function $V(\bx)$ and its gradient $\nabla V(\bx)$. Both these quantities can be obtained as a by-product of solving the SDRE at no additional computational cost since $V(\bx) = \bx^{\top}\Pi(\bx)\bx$ and $\nabla V(\bx) = 2\Pi(\bx)\bx$. 

The sampling datasets are split into \emph{training sets} and \emph{valuation sets}, with a ratio of $80/20$. The goodness of fit in the valuation set, measured by the coefficient of determination $r^2$, guided the choice of the NN's architecture within the FNN family. The minimization of the loss function \eqref{loss} was performed using the quasi-Newton method \texttt{lbfgs}. The parameters to be optimized are the number of hidden layers, the number of neurons per layer, the activation function, and the number of epochs taken into account during the training (we fixed the batches' size to $100$). For $V_{\theta}(\bx)$, we also optimize the hyper-parameters $\mu_V$ and $\mu_{dV}$ weighting the terms in the loss function\eqref{loss1}.
The goodness of fit of the trained models is finally evaluated in the \emph{test set}, a uniform grid of $N=10^4$ points within the state space, where the approximated control is compared with the pointwise computation through the SDRE solution.

\subsection{Test 1: Stabilization for the Cucker-Smale model} 
We test our approach over a high-dimensional, nonlinear and nonlocal control problem related to consensus control of agent-based dynamics: the Cucker-Smale model for consensus control.  We consider $N_a = 20$ agents having states $\bx_i=(y_i,v_i) \in \R^2$, denoting position and velocity respectively, in $\Omega=[-3,3]^{40} \subset \mathbb{R}^{20} \times \mathbb{R}^{20}$ and governed by the dynamics
\begin{align}\label{cs}
    \dot y_i &= v_i\\
    \dot v_i &= \frac{1}{N_a} \sum_{j = 1}^{N_a} \frac{v_j - v_i}{1 + ||y_i - y_j||^2} + u_i
\end{align}
where $i = 1,\ldots,N_a$. Here, the control vector $\bu(t) \in \mathcal{L}^2([0,T];\mathbb{R}^{\times N_a})$ is optimized according to
\begin{equation}
    \min_{\bu(\cdot)} \mathcal{J}(\bx(\cdot)) = \frac{1}{N_a}\int_{0}^T  \sum_{i=1}^{N_a} ||y_i||^2 + ||v_i||^2 + ||u_i||^2 dt
\end{equation}
and it can be written in semilinear form as  
\begin{align*}
	    \begin{bmatrix}
		\dot{\by}\\\dot{\bv}
	\end{bmatrix} &= \begin{bmatrix}
	\mathbb{O}_{N_a} & \mathbb{I}_{N_a}\\\mathbb{O}_{N_a} & \cA_{N_{a}}(\by)\end{bmatrix}
	\begin{bmatrix}
		\by\\\bv
	\end{bmatrix} +\begin{bmatrix}
	\mathbb{O}_{N_a}\\
	\mathbb{I}_{N_a}
\end{bmatrix} \bu\,,\\\\
\big[\cA(\by)\big]_{i,j} &= 
    \begin{cases}
    -\frac{1}{N_a} \sum_{k=1}^{N_a} P(y_i,y_k) &\text{if } i = j\,,\\
    \frac{1}{N_a}  P(y_i,y_j) &\text{otherwise} 
    \end{cases}\\
P(y_i,y_j)&= \frac{1}{1 + ||y_i-y_j||^2}\,,\\
\bQ& = \frac{1}{N_a}\mathbb{I}_{2N_a}, \quad \bR = \mathbb{I}_{N_a} \,,  
\end{align*}
where $\mathbb{O}_{n}$ denotes a matrix of zeros in $\R^{n\times n}$. We train a model for $V_{\theta}$ consisting of a FNN with $3$ hidden  layers with $400$ neurons per layer and activation function $\sigma(x) = \max(0,x)$. The best configuration resulting from hyper-parameter tuning was $(\mu_V,\mu_{dV}) = (0.1, 2)$, where the NN reaches the maximum $r^2$ being trained for 41 epochs, just before overfitting. Finally, applying the trained model to a grid of points in the hypercube $[-3,3]^{40}$, we compute the gradient of the model w.r.t. its input via automatic differentiation, computing the approximate control as in (\ref{u_dV}). The direct feedback model $\bu_{\theta}(\bx) \in \mathbb{R}^{N_a}$ consists of $2$ hidden layers, with $400$ neurons per layer, and activation function $\sigma(x) = tanh(x)$, while being trained for 20 epochs. Goodness of fit for both models is shown in Table \ref{tab2}.

\begin{table}[h]
	\centering 
 \begin{tabular}{c c c}
 	\\
  \emph{predicted variable} & $r^2$ & \emph{MSE}\\ 
  \hline \\
  $V_{\theta}$     & 0.67236    & 0.39829\\ 
  $\nabla V_{\theta}$     & 0.94921    & 0.07906 \\ 
  $\bu_{V}$   & 0.92415    & 56.2208 \\ 
   $\bu_{\theta}$ & 0.96039    & 29.3591 
\end{tabular}
\caption{Goodness of fit for Test 1.}\label{tab2}
\end{table}
\begin{figure}[h]
	\centering
\includegraphics[width=0.22\textwidth]{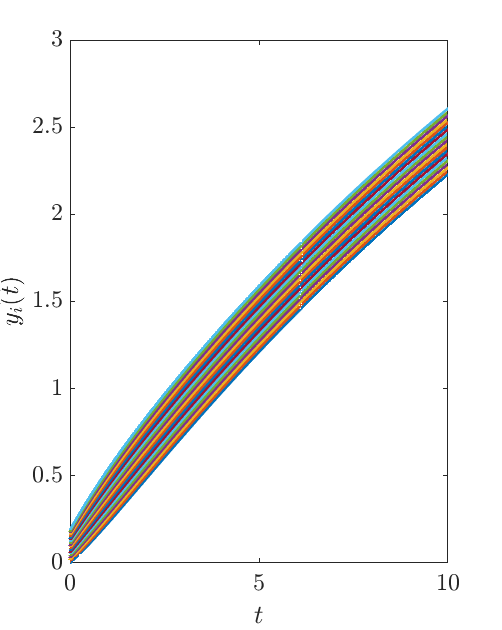}
\includegraphics[width=0.22\textwidth]{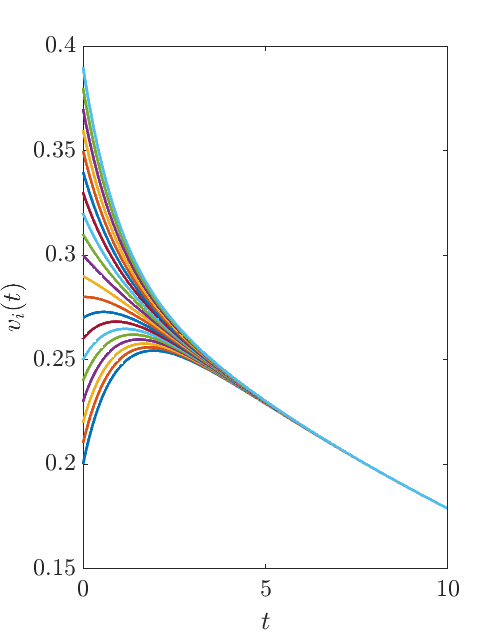}\\
\includegraphics[width=0.22\textwidth]{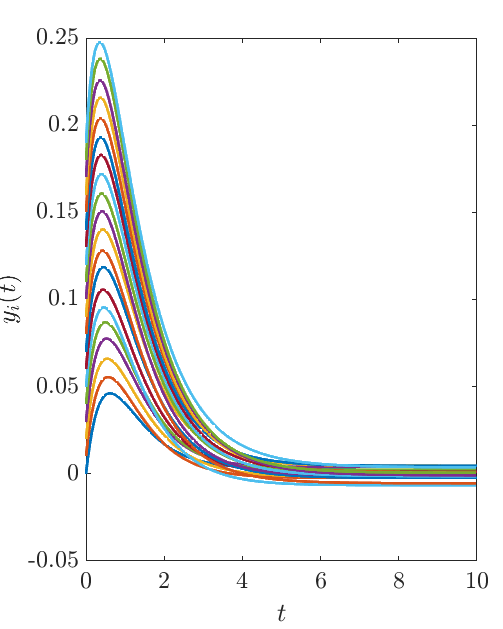}
\includegraphics[width=0.22\textwidth]{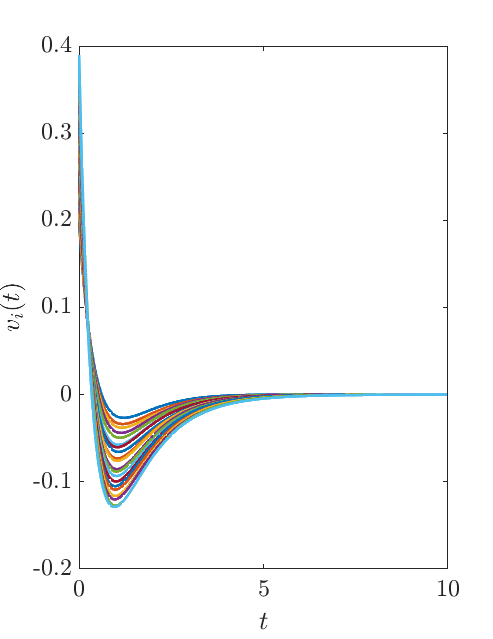}
	\caption{Test 1. Top: uncontrolled positions (left) and velocities (right). Bottom: the controller $\bu_{\theta}(\bx)$ stabilizes the dynamics to the origin.}\label{t2traj}
\end{figure}
\begin{figure}[h]
	\centering
\includegraphics[width=0.24\textwidth]{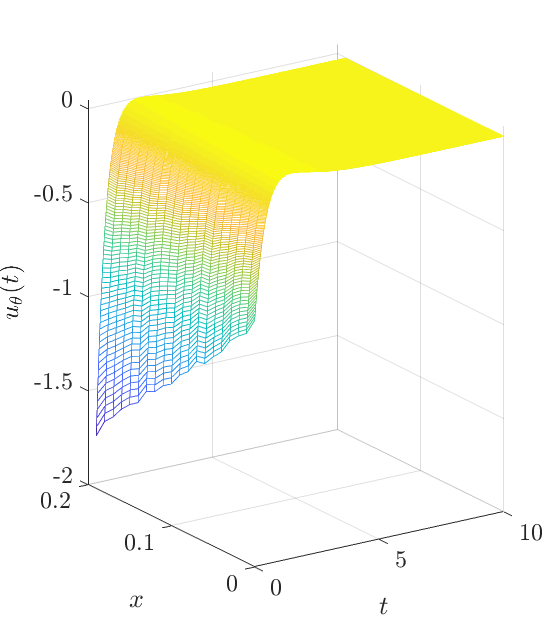}\hfill
\includegraphics[width=0.24\textwidth]{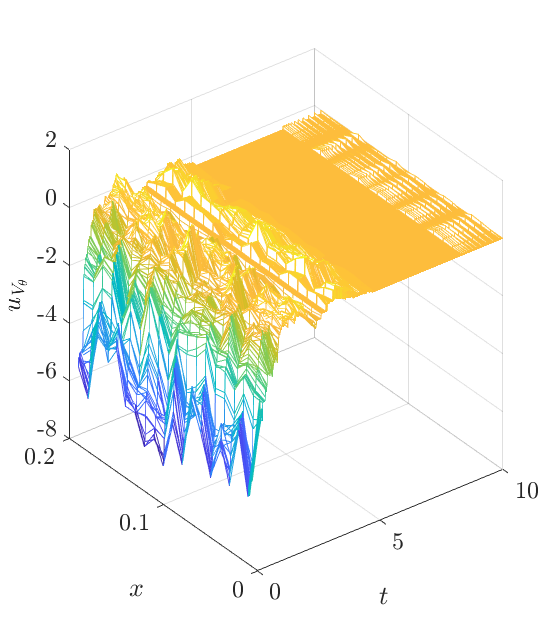}
	\caption{Test 1. Evolution of the control signals $\bu_{\theta}(t)$ (left) and $\bu_{V}(t)$ (right). In this nonlinear test, learning directly a model for $\bu$ leads to better results. We omit the plot for $\bu(t)$ since it is well approximated by $\bu_{\theta}(t)$.}\label{cont2}
\end{figure}
Figure \ref{t2traj}  depicts trajectories for an initial condition $(\by(0),\bv(0))$ that is a vector of equally spaced entries in $[0,0.4]$. For this choice of dimension of the physical space and number of agents, the model $\bu_{\theta}$ performs better than the gradient-augmented $\bu_V$. The differences between both control signals can be observed in Figure \ref{cont2}.

\subsection{Test 2: Feedback control of the Allen-Cahn PDE}
Following a test presented in \cite{tensor4}, we consider the control of the nonlinear Allen-Cahn PDE
\begin{equation}\label{allen_cahn}
    \partial_t x(\xi, t) = 0.1\partial^{2}_{\xi\xi} x+ x(1 + x^2) - \chi_{\omega}(\xi)u(t)
\end{equation}
in $[0,1]\times\mathbb{R}_{+}$ with Neumann boundary conditions, where the scalar control signal $u:[0,+\infty]\to\mathbb{R}$ acts through the indicator function of the interval $\omega = [0.6,0.9]$. Without control action, these dynamics are bistable with $x\equiv \pm 1$ being the stable equilibria. We are interested in minimizing 
\begin{equation}
    \mathcal{J}(u,x) = \int_{0}^{+\infty} ||x(\xi,t)||^2 + 0.1 u^2(t) dt\,,
\end{equation}
thus stabilizing the dynamics towards the equilibrium $x = 0$. 
The PDE (\ref{allen_cahn}) is discretized in space via finite differences with $N=51$ degrees of freedom, leading to a system of non-linear ODEs
\begin{equation}\label{disc_ac}
    \dot\bx = \bA\bx + \bx \odot (1 + \bx \odot \bx) + \bB u(t),
\end{equation}
where $\bx(t) = (x(\xi_1,t),...,x(\xi_N,t))$ is the discrete state, $\odot$ denotes the Hadamard product and $\bA,\bB$ correspond to a discretization of the Laplace operator and the indicator function $\chi_{\omega}(\xi)$ over a uniform grid $\{\xi_i\}_{i=1}^N$.

We consider a dataset $\{\bx^{(i)},V(\bx^{(i)}),\nabla V(\bx^{(i)})\}_{i=1}^{N_s}$ with $N_s = 1000$, where the states have been sampled from $[-2,2]^{51}$. We train a model for $V_{\theta}$ with $3$ hidden layers, $500$ neurons per layer, and activation function $\sigma(x) = \max(0,x)$. The best configuration of hyper-parameters is found to be $(\mu_V,\mu_{dV}) = (0.9, 7)$, with the NN being trained for 71 epochs. Finally, we test the trained model $V_{\theta}$ in a test grid of points in $[-2,2]^{51}$. For the model $\bu_{\theta}$, the architecture is built with $4$ hidden layers, with $500$ neurons per layer, and activation function $\sigma(x) = \max(0,x)$. In this example, also the output layer for $\bu_{\theta}$ is made only of a single neuron, since the feedback law is a scalar in $\mathbb{R}$. The model was trained for 50 epochs.
\begin{table}[h]
\begin{center}
 \begin{tabular}{c c c} 
  \emph{predicted variable} & $r^2$ & \emph{MSE}\\ 
  \hline\\
  $V_{\theta}$     & 0.81681       & 0.00024628\\ 
  $\nabla V_{\theta}$    & 0.87114       & 0.00025852\\ 
  $\bu_{V}$ & 0.91976      & 0.012181\\ 
  $\bu_{\theta}$ & 0.85443   & 0.022098 
\end{tabular}
\end{center}
\caption{Goodness of fit for Test 2.}\label{tabt3}
\end{table}

\begin{figure}[!ht]
	\centering
	\includegraphics[width=8cm]{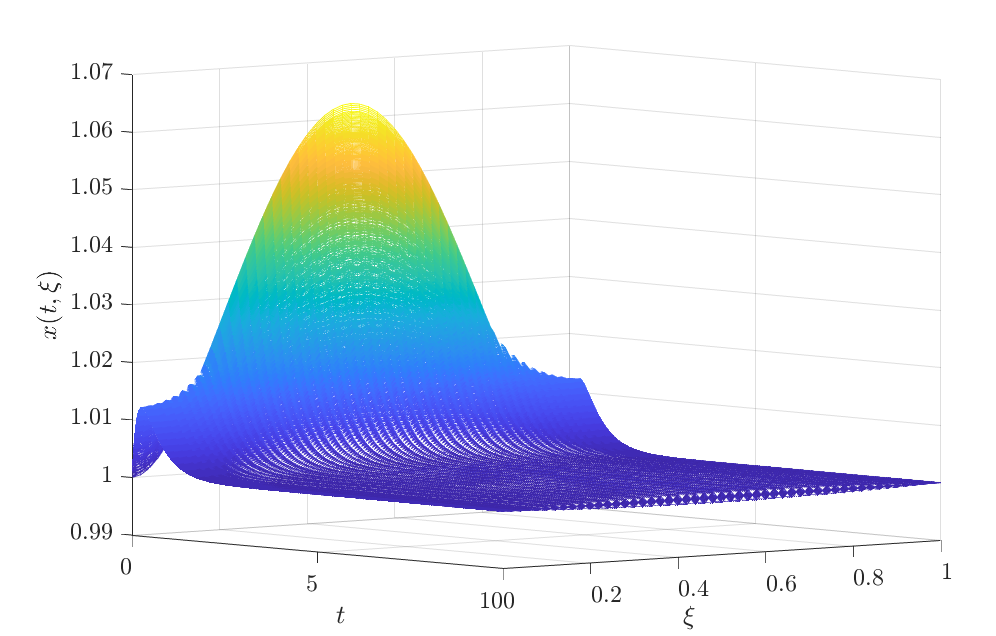}\\
	\includegraphics[width=8cm]{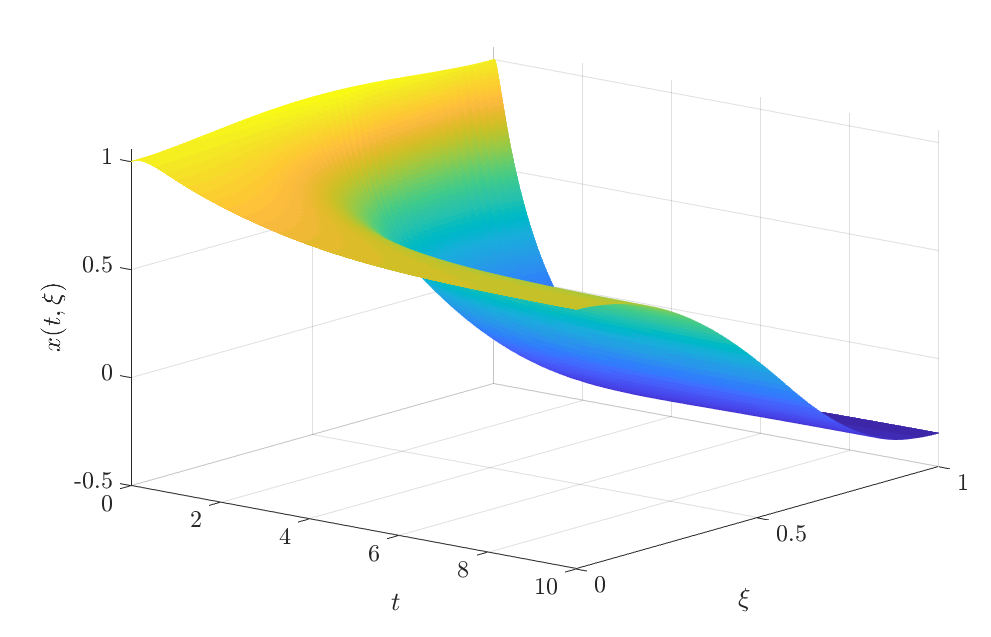}\\
	\includegraphics[width=8cm]{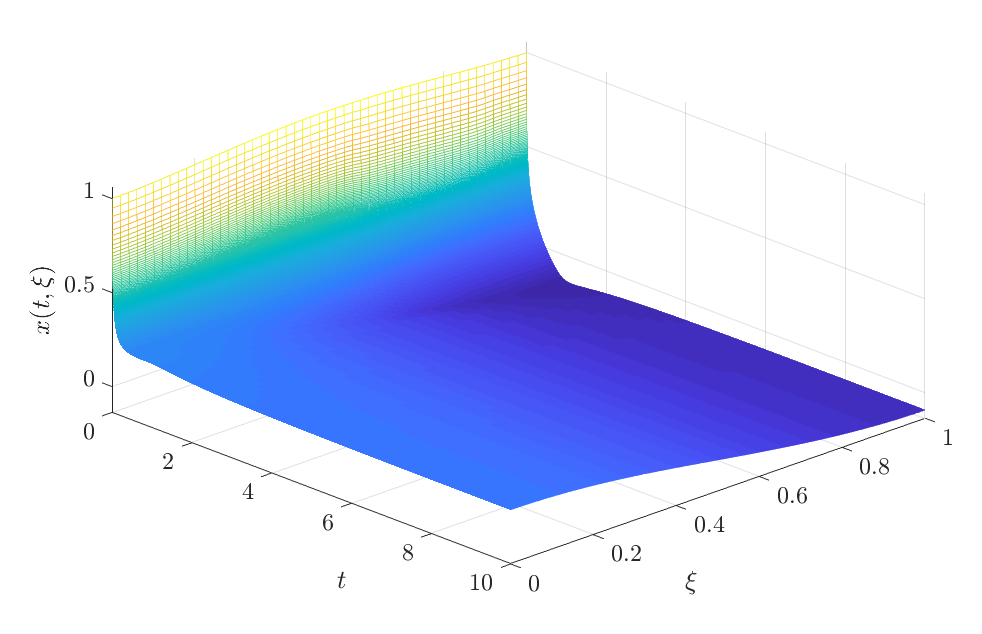}\\
	\caption{Test 2. Top: uncontrolled state, converges to $x=1$. Middle: controlled state with linearized feedback $\bu_0(\bx)=-\bK_0\bx$ around the origin, fails to stabilize. Bottom: the nonlinear feedback $\bu_V$ stabilizes the system towards the origin.}\label{trajt3}
\end{figure}

\begin{figure}[h]
	\centering
	\includegraphics[width=0.4\textwidth]{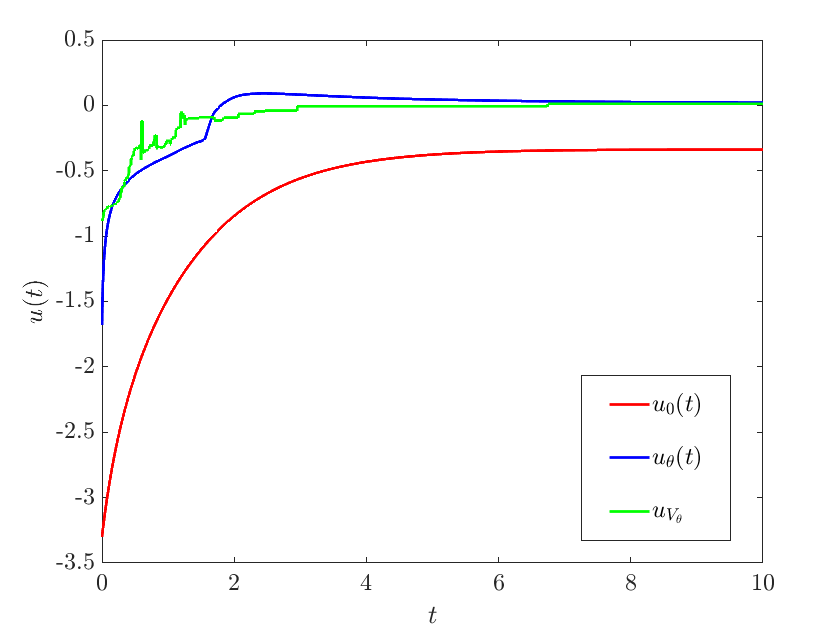}
	\caption{Test 2. Different control signals: $\bu_0$, $\bu_{\theta}$ and $\bu_V$. In this test the feedback law trained with gradient-augmented data $\bu_V$ outperforms $\bu_{\theta}$.}\label{cont3}
\end{figure}

In Figs \ref{trajt3} and \ref{cont3} we compare the trajectories resulting from the integration of the discretized dynamics (\ref{disc_ac}) with $t\in[0,10]$, for an initial condition $x(\xi,0)=1+(1-\xi)\xi$, and different feedback laws: the constant zero function, the feedback obtained considering the linear control operator $\bK_0$, the control resulting from the gradient-augmented approximation $V_{\theta}$, and the one given by $\bu_{\theta}$. In this high-dimensional local problem, the approximation done through the gradient-augmented model $V_{\theta}$ happens to outperform $\bu_{\theta}$ in terms of goodness of fit. On the other hand, observing the different closed-loop evolutions and control signals, we can see how both approximate feedback laws succeed in stabilizing the trajectories near $x = 0$, while the uncontrolled system is stable in $x = 1$ and the $\bu_0$ results in a system which, for $t=10$ has not yet approached the equilibrium.

\section{Conclusions}
We have presented a novel computational method for the approximation of stabilizing feedback laws in nonlinear dynamics based on  a supervised learning approach. The training data originates from the pointwise solution of the State-Dependent Riccati Equation. We have studied the approximation of the feedback control through feedforward neural networks, and analysed different choices of architectures and loss functions for training. We have provided computational evidence that for high-dimensional nonlinear problems, the SDRE feedback law can be effectively approximated through FNNs, thus removing the stringent requirement of a fast ARE solver for real-time closed-loop control. We proposed two alternatives for learning a model for the control. We observe that for genuinely nonlinear control problems, such as agent-based dynamics, better results are achieved by learning directly the feedback $\bu(\bx)$ from the SDRE solves. However, for problems where a linear structure is more prominent, such as in the control of semilinear parabolic PDEs, learning a model for a local approximation of the value function $V(\bx)$ and computing the control from its gradient is an accurate and more efficient alternative.

\addtolength{\textheight}{-12cm}   




\bibliographystyle{IEEEtran}
\bibliography{IEEEabrv,polynomshort}

\end{document}